\documentclass[runningheads]{cl2emult}

\usepackage{makeidx}  
\usepackage{graphicx} 
\usepackage{subeqnar} 
\usepackage{multicol} 
\usepackage{cropmark} 
\usepackage{math}     

\usepackage{amssymb}
\usepackage{amsmath}   
\newtheorem{thm}{Theorem}
\newtheorem{lem}{Lemma}

\newtheorem{cor}{Corollary}
\newtheorem{conj}{Conjecture}
\newtheorem{ex}{Example}

\begin{document}
\title*{Multiplicities of Singular Points in Schubert Varieties
of Grassmannians}
\titlerunning{Singular Points in Schubert Varieties}
%
%
%
%
%
\author{Victor Kreiman \and V. Lakshmibai}
\institute{Northeastern University, Boston, MA 02115}
\authorrunning{Victor Kreiman \and V. Lakshmibai}
\maketitle

%
%

\begin{abstract}
\index{abstract} We give a closed-form formula for the Hilbert
function of the tangent cone at the identity of a Schubert variety
$X$ in the Grassmannian in both group theoretic and combinatorial
terms. We also give a formula for the multiplicity of $X$ at the
identity, and a Gr\"obner basis for the ideal defining $X \cap
O^-$ as a closed subvariety of $O^-$, where $O^-$ is the opposite
cell in the Grassmannian. We give conjectures for the Hilbert
function and multiplicity at points other than the identity.
\end{abstract}

\section{Introduction}

The first formulas for the multiplicities of singular points on
Schubert varieties in Grassmannians appeared in Abhyankar's
results \cite{ab} on the Hilbert series of determinantal varieties
(recall that a determinantal variety gets identified with the
opposite cell in a suitable Schubert variety in a suitable
Grassmannian). Herzog-Trung \cite{herzog} generalized these
formulas to give determinantal formulas for the multiplicities at
the identity of all Schubert varieties in Grassmannians. Using
standard monomial theory, Lakshmibai-Weyman \cite{l-w} obtained a
recursive formula for the multiplicities of all points in Schubert
varieties in a minuscule $G/P$; Rosenthal-Zelevinsky
\cite{rose-zel} used this result to obtain a closed-form
determinantal formula for multiplicities of all points in
Grassmannians.

\section{Summary of Results}\label{summary}
Let $K$ be the base field, which we assume to be algebraically
closed, of arbitrary characteristic. Let $G$ be $SL_n(K)$, $T$ the
subgroup of diagonal matrices in $G$, and $B$ the subgroup of
upper diagonal matrices in $G$.  Let $R$ be the root system of $G$
relative to $T$, and $R^+$ the set of positive roots relative to
$B$.  Let $W$ be the Weyl group of $G$.  Note that $W = S_n$, the
group of permutations of the set of $n$ elements.  Let $P_d$ be
the maximal parabolic subgroup
\begin{equation*}
P_d=\left\{ A\in G\biggm| A=
\begin{pmatrix}
*&*\\ 0_{(n-d)\times d}&*
\end{pmatrix}\right\}.
\end{equation*}
Let $R_{P_d}$, $R_{P_d}^+$, and $W_{P_d}$ denote respectively the
root system, set of positive roots, and Weyl group of ${P_d}$. The
quotient $W/W_{P_d}$, with the Bruhat order, is a distributive
lattice. The map $\alpha \mapsto s_\alpha W_{P_d}$ taking a
positive root to its corresponding reflection, embeds $R^+
\!\setminus\! R_{P_d}^+$ in $W/W_{P_d}$. We shall also denote the
image by $R^+ \!\setminus\! R_{P_d}^+$. It is a sublattice of
$W/W_{P_d}$.

A {\it multiset} is similar to a set, but with repetitions of
entries allowed.  Define the {\it cardinality} of a multiset $S$,
denoted by $|S|$, to be the number of elements in $S$, including
repetitions. Define a {\it uniset} to be a multiset which has no
repetitions. If $\mathcal{S}$ is a set, define
$\mathcal{S}^{\ast}$ to be the collection of all multisets which
are made up of elements of $\mathcal{S}$.

A {\it chain of commuting reflections} in $W / W_{P_d}$ is a
nonempty set of pairwise-commuting reflections $\{s_{\alpha_1},
\ldots, s_{\alpha_t}\}$, $\alpha_i \in R^+ \!\setminus\!
R_{P_d}^+$, such that $s_{\alpha_1} > \cdots > s_{\alpha_t}$; we
refer to $t$ as the {\it length} of the chain. For a multiset $S
\in (R^+ \!\setminus\! R_{P_d}^+)^{\ast}$, define the {\it
chainlength} of $S$ to be the maximum length of a chain of
commuting reflections in $S$.

Fix $w \in W/W_{P_d}$.  Define $S_w$ to be the multisets $S$ of
$(R^+ \!\setminus\! R_{P_d}^+)^{\ast}$, such that the product of
every chain of commuting reflections in $S$ is less than or equal
to $w$; similarly, define $S_w'$ to be the unisets of $(R^+
\!\setminus\! R_{P_d}^+)^{\ast}$ having the same property.  For
$m$ a positive integer, define $$S_w(m) = \{ S\in S_w \, : \,
|S|=m\}$$ $$S_w'(m) = \{ S\in S_w' \, : \, |S|=m\}.$$

We can now state our two main results.  First, letting $X(w)$
denote the Schubert variety of $G/{P_d}$ corresponding to $w \in
W/W_{P_d}$, the Hilbert function of the tangent cone to $X(w)$ at
the identity is given by

\begin{thm}\label{main1}
$h_{\textrm{TC}_{id}X(w)}(m) = |S_w(m)|, m \in \mathbb{N}$.
\end{thm}
Second, letting $M$ denote the maximum cardinality of any element
of $S_w'$, the multiplicity at the identity is given by

\begin{thm}\label{main2}
$\textrm{mult}_{id}X(w) = |\{S \in S_w' \, : \, |S|=M\}|$.
\end{thm}

\section{Preliminaries}

\subsection{Multiplicity of an Algebraic Variety at a Point}

Let $B$ be a graded, affine $K$-algebra such that $B_1$ generates
$B$ (as a $K$-algebra). Let $X = \textrm{Proj}(B)$. The function
$h_B(m)$ (or $h_X(m)$) $= \textrm{dim}_KB_m$, $m \in \mathbb{Z}$
is called the {\it Hilbert function} of $B \,(\textrm{or } X)$.
There exists a polynomial $P_B(x)$ (or $P_X(x)$) $\in
\mathbb{Q}[x]$, called the {\it Hilbert polynomial} of $B\,
(\textrm{or } X)$, such that $f_B(m) = P_B(m)$ for $m \gg 0$. Let
$r$ denote the degree of $P_B(x)$. Then $r = \textrm{dim}(X)$, and
the leading coefficient of $P_B(x)$ is of the form $c_B/r!$, where
$c_B \in \mathbb{N}$. The integer $c_B$ is called the {\it degree
of $X$}, and denoted $\textrm{deg}(X)$. In the sequel we shall
also denote $\textrm{deg}(X)$ by $\textrm{deg}(B)$.

Let $X$ be an algebraic variety, and let $P \in X$.  Let $A =
\mathcal{O}_{X, P}$ be the stalk at $P$ and $\mathfrak{m}$ the
unique maximal ideal of the local ring $A$. Then the {\it tangent
cone} to $X$ at $P$, denoted $\textrm{TC}_P(X)$, is defined to be
$\textrm{Spec}(\textrm{gr}(A,\mathfrak{m}))$, where
$\textrm{gr}(A,\mathfrak{m}) =
\oplus_{j=0}^{\infty}\mathfrak{m}^j/\mathfrak{m}^{j+1}$. The {\it
multiplicity} of $X$ at $P$, denoted $\textrm{mult}_P(X)$, is
defined to be
$\textrm{deg}(\textrm{Proj}(\textrm{gr}(A,\mathfrak{m})))$. If $X
\subset K^n$ is an affine closed subvariety, and $m_P \subset
K[X]$ is the maximal ideal corresponding to $P \in X$, then
$\textrm{gr}(K[X], m_P) = \textrm{gr}(A,\mathfrak{m})$.

\subsection{Monomial Orders, Gr\"obner Bases, and Flat Deformations}

Let $A$ be the polynomial ring $K[x_1, \cdots, x_n]$.  A {\it
monomial order} $\succ$ on the set of monomials in $A$ is a total
order such that given monomials $m$, $m_1$, $m_2$, $m \neq 1$,
$m_1 \succ m_2$, we have $mm_1 \succ m_1$ and $mm_1 \succ mm_2$.
The largest monomial (with respect to $\succ$) present in a
polynomial $f \in A$ is called the {\it initial term} of $f$, and
is denoted by $\textrm{in}(f)$.

The {\it lexicographic order} is a total order defined in the
following manner.  Assume the variables $x_1, \ldots, x_n$ are
ordered by $x_n > \cdots > x_1$. A monomial $m$ of degree $r$ in
the polynomial ring $A$ will be written in the form $m = x_{i_1}
\cdots x_{i_r}$, with $n \ge i_1 \ge \cdots \ge i_r \ge 1$.  Then
$x_{i_1} \cdots x_{i_r} \succ x_{j_1} \cdots x_{j_s}$ in the
lexicographic order if and only if either $r > s$, or $r = s$ and
there exists an $l < r$ such that $i_1 = j_1, \ldots, i_l = j_l,
i_{l+1} > j_{l+1}$.  It is easy to check that the lexicographic
order is a monomial order.

Given an ideal $I \subset A$, denote by $\textrm{in}(I)$ the ideal
generated by the initial terms of the elements in $I$.  A finite
set $\mathcal{G} \subset I$ is called a {\it Gr\"obner basis} of
$I$ (with respect to the monomial order $\succ$), if
$\textrm{in}(I)$ is generated by the initial terms of the elements
of $\mathcal{G}$.

\subsubsection{Flat Deformations: }\label{flatdef}

Given a monomial order and an ideal $I \subset A$, there exists a
flat family over $\textrm{Spec}(K[t])$ whose special fiber $(t =
0)$ is $\textrm{Spec}(A/\textrm{in}(I))$ and whose generic fiber
($t$ invertible) is $\textrm{Spec}(A/I \otimes K[t, t^{-1}])$.
Further, if $I$ is homogeneous, then the special fiber and generic
fiber have the same Hilbert function (see \cite{eis} for details).

\subsection{Grassmannian and Schubert Varieties}

\subsubsection{The Pl\"ucker Embedding: }

Let $d$ be such that $1 \le d < n$. The {\it Grassmannian}
$G_{d,n}$ is the set of all $d-$dimensional subspaces $U \subset
K^n$.  Let $U$ be an element of $G_{d,n}$ and $\{ a_1, \ldots a_d
\}$ a basis of $U$, where each $a_j$ is a vector of the form
\vspace{-.05in}
\begin{displaymath}
a_j = \left(\begin{array}{c} a_{1j} \\ a_{2j} \\ \vdots \\ a_{nj}
\end{array}\right),\, \textrm{with } a_{ij} \in K.
\end{displaymath}
Thus, the basis $\{a_1, \cdots, a_d\}$ gives rise to an $n \times
d$ matrix $A = (a_{ij})$ of rank $d$, whose columns are the
vectors $a_1, \cdots, a_d$.

We have a canonical embedding
\begin{equation*}
p:G_{d,n}\hookrightarrow \mathbb{P}(\wedge^dK^n) \ , \ U \mapsto
[a_1 \wedge \cdots \wedge a_d]
\end{equation*}
called the {\it Pl\" ucker embedding}.  Let
\begin{equation*}
I_{d,n}=\{\underline{i}=(i_1,\dots,i_d) \in \mathbb{N}^d \, : \,
1\le i_1<\dots<i_d\le n\}\ .
\end{equation*}
Then the projective coordinates ({\it Pl\"ucker coordinates}) of
points in $\mathbb{P}(\wedge^d K^n)$ may be indexed by $I_{d,n}$;
for $\underline{i} \in I_{d,n}$, we shall denote the
$\underline{i}$-th component of $p$ by $p_{\underline{i}}$, or
$p_{i_1, \cdots, i_d}$.  If a point $U$ in $G_{d, n}$ is
represented by the $n \times d$ matrix $A$ as above, then $p_{i_1,
\cdots, i_d}(U) = \textrm{det}(A_{i_1, \ldots, i_d})$, where
$A_{i_1, \ldots, i_d}$ denotes the $d \times d$ submatrix whose
rows are the rows of $A$ with indices $i_1, \ldots, i_d$, in this
order.

\subsubsection{Identification of $G/P_d$ with $G_{d,n}$: }

Let $G$, $T$, $B$, and $P_d$ be as in Section~\ref{summary}. Let
$\{e_1, \ldots, e_n\}$ be the standard basis for $K^n$.  For the
natural action of $G$ on $\mathbb{P}(\wedge^d K^n)$, the isotropy
group at $[e_1\wedge \cdots \wedge e_d]$ is $P_d$, while the orbit
through $[e_1\wedge \cdots \wedge e_d]$ is $G_{d,n}$.  Thus we
obtain an identification of $G / P_d$ with $G_{d, n}$.  We also
note that $W/W_{P_d} \,( \, =S_n / (S_d \times S_{n-d}))$ may be
identified with $I_{d,n}$.

\subsubsection{Schubert Varieties: }

For the action of $G$ on $G_{d, n}$, the $T$-fixed points are
precisely $\{ [e_{\underline{i}}\,], \underline{i} \in I_{d, n}
\}$, where $e_{\underline{i}} = e_{i_1} \wedge \cdots \wedge
e_{i_d}$ . The {\it Schubert variety} $X_{\underline{i}}$
associated to $\underline{i}$ is the Zariski closure of the
$B$-orbit $B[e_{\underline{i}} \,]$ with the canonical reduced
scheme structure.

We have a bijection between $\{ \textrm{Schubert varieties in
}G_{d, n} \}$ and $I_{d, n}$.  The partial order on Schubert
varieties given by inclusion induces a partial order (called the
{\it Bruhat order}) on $I_{d, n} \, ( \, = W / W_{P_d})$; namely,
given $\underline{i}= (i_{1},\dots,i_{d})$, $\underline{j}=
(j_{1},\dots ,j_{d}) \in I_{d,n}$,
\begin{equation*}
\underline{i}\ge \underline{j} \iff i_t\ge j_t, \text{for all } 1
\leq t\leq d.
\end{equation*}

We note the following facts for Schubert varieties in the
Grassmannian (see \cite{g-l2} or \cite{musili} for example):

\begin{itemize}
\item \textbf{Bruhat Decomposition}:
$\displaystyle X_{\underline{i}}=
\bigcup_{\underline{j}\,\le\,\,\underline{i}}B[e_{\underline{j}}]$.
\item \textbf{Dimension}: $\displaystyle \dim X_{\underline{i}}=\!\!\sum_{1\le t\le
d}\!i_t-t$.
\item \textbf{Vanishing Property of a Pl\"ucker Coordinate}:
$$p_{\underline{j}}\bigr|_{X_{\underline{i}}}\ne
0\iff\underline{i}\ge\underline{j}.$$
\end{itemize}

\subsubsection{Standard Monomials: }

A monomial $f = p_{\theta_1}\!\! \cdots p_{\theta_t}$, $\theta_i
\in W / W_{P_d}$ is said to be {\it standard} if \vspace{-.05in}
\begin{equation}\label{standard}
\theta_1 \geq \cdots  \geq \theta_t.
\end{equation}
Such a monomial is said to be {\it standard on the Schubert
variety $X(\theta)$}, if in addition to (\ref{standard}), we have
$\theta \geq \theta_1$.

Let $w \in W / W_{P_d}$.  Let $R(w) = K[X(w)]$, the homogeneous
coordinate ring for $X(w)$, for the Pl\"ucker embedding.  Recall
the following two results from standard monomial theory (cf.
\cite{g-l2}).
\begin{thm}\label{rwbasis}
The set of standard monomials on $X(w)$ of degree $m$ is a basis
for $R(w)_m$.
\end{thm}

\begin{thm}\label{iwgen}
For $w \in W / W_{P_d}$, let $I_w$ be the ideal in $K[G_{d, n}]$
generated by $\{ p_{\theta}, \theta \nleq w \}$.  Then $R(w) =
K[G_{d,n}]/I_w$.
\end{thm}

\subsubsection{The Opposite Big Cell $O^-$: }

Let $U^-$ denote the unipotent lower triangular matrices of
$G=SL_n(K)$. Under the canonical projection $G \to G/P_d \, , \, g
\mapsto gP_d \, ( \, = g[e_{id}])$, $U^-$ maps isomorphically onto
its image $U^-[e_{id}]$.  The set $U^-[e_{id}]$ is called the {\it
opposite big cell} in $G_{d, n}$, and is denoted by $O^-$. Thus,
$O^-$ may be identified with
\begin{equation}\label{ominus}
\left\{
\begin{pmatrix} &\text{\rm Id}_{d\times d}\\ x_{d+1\, 1} &\dots &x_{d+1\, d}\\
\vdots &&\vdots\\ x_{n\, 1}&\dots&x_{n\,d}
\end{pmatrix},\quad x_{ij}\in K,\quad d+1\le i\le n, 1\le j\le d
\right\}.
\end{equation}
Thus we see that $O^-$ is an affine space of dimension $(n-d)
\times d$, with $id$ as the origin; further $K[O^-]$ can be
identified with the polynomial algebra $K[x_{-\beta},\beta\in R^+
\setminus R_{P_d}^+]$.  To be very precise, denoting the elements
of $R$ as in \cite{bou}, we have $R^+ \setminus R_{P_d}^+
=\{\epsilon_j-\epsilon_i$,\, $d+1\le i\le n, \, 1\le j\le d\}$;
given $\beta\in R^+ \setminus R_{P_d}^+$, say
$\beta=\epsilon_j-\epsilon_i$, we identify $x_{-\beta}$ with
$x_{ij}$.  We denote by $s_{(i,j)} \, (\textrm{or } s_{(j,i)})$
the reflection corresponding to $\beta$, namely, the transposition
switching $i$ and $j$.

\subsubsection{Evaluation of Pl\"ucker Coordinates on $O^-$: }

Let $\underline{j} \in I_{d, n}$.  We shall denote the Pl\"ucker
coordinate $p_{\underline{j}}|_{O^-}$ by $f_{\underline{j}}$.  Let
us denote a typical element $A \in O^-$ by $\begin{pmatrix}
\text{\rm Id}_{d\times d}\\ X \end{pmatrix}$. Then
$f_{\underline{j}}$ is simply a minor of $X$ as follows. Let
$\underline{j} = (j_1, \ldots, j_d)$, and let $j_r$ be the largest
entry $ \leq d$. Let $\{k_1, \ldots, k_{d-r}\}$ be the complement
of $\{j_1, \ldots, j_r\}$ in $\{ 1, \ldots, d \}$. Then this minor
of $X$ is given by column indices $k_1, \ldots, k_{d-r}$ and row
indices $j_{r+1}, \ldots, j_d$ (here the rows of $X$ are indexed
as $d+1, \ldots, n$).

Conversely, given a minor of $X$, say, with column indices $b_1,
\ldots, b_s$, and row indices $i_{d-s+1}, \ldots, i_d$, then that
minor is the evaluation of $f_{\underline{j}}$ at $X$, where
$\underline{j} = (j_1, \ldots, j_d)$ may be described as follows:
$\{j_1, \ldots, j_{d-s}\}$ is the complement of $\{b_1, \ldots,
b_s \}$ in $\{1, \ldots, d \}$, and $j_{d-s+1}, \ldots, j_d$ are
simply the row indices (again, the rows of $X$ are indexed as
$d+1, \ldots, n$).

Note that if $\underline{j} = (1, \ldots, d)$, then
$p_{\underline{j}}$ evaluated at $X$ is 1. In the above
discussion, therefore, we must consider the element 1 (in
$K[O^-]$) as the minor of $X$ with row indices (and column
indices) given by the empty set.

\begin{ex} Consider $G_{2,4}$. Then
$$ O^- = \left\{
\begin{pmatrix} 1&0\\ 0&1\\ x_{31}&x_{32}\\ x_{41}&x_{42}
\end{pmatrix},\; x_{ij}\in K
\right\}. $$ On $O^-$, we have $p_{12}=1$, $p_{13}=x_{32}$,
$p_{14}=x_{42}$, $p_{23}=x_{31}$, $p_{24}=x_{41}$,
$p_{34}=x_{31}x_{42}-x_{41}x_{32}$.
\end{ex}
Note that each of the Pl\"ucker coordinates is homogeneous in the
local coordinates $x_{ij}$.

\section{The Hilbert Function of $\textrm{TC}_{id}X(w)$}

In view of the Bruhat decomposition, in order to determine the multiplicity
at a singular point $x$, it is enough to determine the
multiplicity of the $T$-fixed point in the $B$ orbit $Bx$.  In
this section, we shall discuss the behavior at a particular
$T$-fixed point, namely the identity.

\subsection{The Variety Y(w)}
We define $Y(w) \subset G_{d,n}$ to be $X(w) \cap O^-$. Since
$Y(w) \subset X(w)$ is open dense, and $id \in Y(w)$, we have that
$\textrm{TC}_{id}Y(w) = \textrm{TC}_{id}X(w)$. As a consequence of
Theorem \ref{iwgen}, $Y(w) \subset O^-$ is defined as an algebraic
subvariety by the homogeneous polynomials $f_{\theta}, \theta
\nleq w$; further, $id \in O^-$ corresponds to the origin. Thus we
have that $\textrm{gr}(K[Y(w)],m_{id}) = K[Y(w)]$. Hence,
\begin{eqnarray}\label{tc}
\textrm{TC}_{id}X(w) = \textrm{TC}_{id}Y(w) &=&
\textrm{Spec}(\textrm{gr}(K[Y(w)],m_{id})){} \nonumber\\ & & {}
=\textrm{Spec}(K[Y(w)]) = Y(w).
\end{eqnarray}

\subsection{Monomials and Multisets}\label{monmult}

For a monomial $p = x_{\alpha_{i_1}} \cdots x_{\alpha_{i_m}}\!\!\!
\in K[O^-]$, define {\it $\textrm{Multisupp}(p)$} to be the
multiset $\{ \alpha_{i_1}, \ldots \alpha_{i_m} \}$.  It follows
immediately from the definition that Multisupp gives a bijection
between the monomials of $K[O^-]$ and the multisets of $(R^+
\!\setminus\! R_{P_d}^+)^{\ast}$, pairing the square-free
monomials with the unisets. Let $w \in W/W_{P_d}$.  We call a
monomial {\it w-good} if it maps under Multisupp to an element of
$S_w$.  Note that the w-good square-free monomials are precisely
those which map to $S_w'$.

Define a monomial order $\succ$ on $K[O^-]$ in the following
manner.  We say $x_{i,j} > x_{i'\!,j'}$ if $i > i'$, or if $i =
i'$ and $j < j'$.  Note that this extends the partial order
$x_{\alpha}
> x_{\beta} \iff s_{\alpha}
> s_{\beta}$ (in the Bruhat order). The monomials are then
ordered using the lexicographic order.

Define the monomial ideal $J_w \subset K[O^-]$ to be the ideal
gen\-er\-ated by $\{ \textrm{in}f_{\theta}, \theta \nleq w\}$, and
let $A_w = K[O^-]/J_w$. With our or\-der\-ing,
$\textrm{Multisupp}(\textrm{in}f_{\theta})$ is a commuting chain
of reflections whose product is $\theta$. Thus the non w-good
monomials form a vector space basis for $J_w$, and therefore the
w-good monomials form a basis for $A_w$.

\subsection{Sketch of Proof of Theorems 1 and 2}

In view of (\ref{tc}) and the above discussion, Theorem
\ref{main1} follows immediately from
\begin{lem}\label{main}
$h_{K[Y(w)]}(m) = h_{A_w}(m) \, , m \in \mathbb{N}$.
\end{lem}
Theorem \ref{main2} is also a consequence.  Indeed,
\begin{equation*}
\textrm{mult}_{id}X(w) = \textrm{deg}(K[\textrm{TC}_{id}X(w)]) =
\textrm{deg}(K[Y(w)]) = \textrm{deg}(A_w).
\end{equation*}
Since $A_w$ is an affine quotient of an ideal generated by
square-free monomials, letting $M$ be the maximum degree of a
square-free monomial in $A_w$, we have (cf. \cite{Cox})
\begin{eqnarray*}
&\,&\textrm{deg}(A_w) = |\{p \in A_w \, : \,   p \textrm{ is a
square-free monomial} \textrm{ and deg}(p)=M\}| \\&\,&\quad = |\{p
\in K[O^-] \, : \, p \textrm{ is a square-free w-good monomial}
\textrm{ and deg}(p)=M\}| \\&\,&\quad = |\{S \in S_w' \, : \,
|S|=M \}|,
\end{eqnarray*}
yielding Theorem \ref{main2}.

The proof of Lemma \ref{main} relies on an inductive argument
which shows directly that both functions agree for all positive
integers $m$. Note that $K[Y(w)] = K[X(w)]_{(p_{id})}$.  Thus, as
a consequence of Theorem \ref{rwbasis}, $K[Y(w)]$ has a basis
consisting of monomials of the form $f_{\theta_1} \cdots
f_{\theta_t}, w \ge \theta_1 \ge \cdots \ge \theta_t$. If
$SM_w(m)$ denotes the basis elements of degree $m$, then
$h_{Y(w)}(m) = |SM_w(m)|$. Letting $d = d_w$ be the degree of $w$
(see section \ref{combin} below for definition), as a consequence
of standard monomial theory we have
\begin{equation}\label{smred}
SM_w(m+d) = SM_w(m)\,\dot{\cup}\,SM_H(m+d)
\end{equation}
where $SM_H(m+d)= \bigcup_{w_i} SM_{w_i}(m+d)$, the union being
taken over the divisors $X(w_i)$ of $X(w)$ (cf. \cite{l-w}).

We have that $|SM_H(m+d)| = |\bigcup_{w_i} SM_{w_i}(m+d)|$ can be
set-theoretically written as the integral linear combination of
terms of the form $|SM_{w_i}(m+d)|$ and terms of the form
$|SM_{w_j}(m+d) \cap \cdots \cap SM_{w_k}(m+d)|$. Further, it can
be shown that
\begin{equation*} SM_{w_j}(m+d) \cap \cdots \cap SM_{w_k}(m+d) =
SM_{\theta}(m+d),
\end{equation*}
where $\theta$ is given by $X(\theta) =  X(w_j) \cap \cdots \cap
X(w_k)$.  (Note that $I_{d, n}$ being a distributive lattice
implies that for $\tau, \phi \in I_{d, n}$, $X(\tau) \cap X(\phi)$
is irreducible.) Thus, \vspace{-.13in}
\begin{equation}\label{smh}
|SM_H(m+d)| = \!\! \sum_{w' < w} \! a_{w'}|SM_{w'}(m+d)|, \textrm{
for some } a_{w'} \in \mathbb{Z}.
\end{equation}
Taking cardinalities of both sides of (\ref{smred}), we obtain
\begin{equation*}
h_{K[Y(w)]}(m+d) = h_{K[Y(w)]}(m) + \!\! \sum_{w' < w} \!a_{w'
\!}h_{K[Y(w')]}(m+d).
\end{equation*}
Equivalently, $h_{K[Y(w)]}$ satisfies the difference equation
\begin{equation}\label{diffeq}
\phi(w,m+d)=\phi(w,m)+ \!\! \sum_{w' < w} \!a_{w'}\phi(w',m+d).
\end{equation}
To prove Lemma \ref{main}, it suffices to show that $h_{A_w}(m)$
satisfies (\ref{diffeq}) for all $m \in \mathbb{Z}_{\geq\, 0}$,
since it is a straightforward verification that $h_{K[Y(w)]}(m)$
and $h_{A_w}(m)$ have the same initial conditions.

As stated earlier, $K[A_w]$ has as basis the w-good monomials of
$K[O^-]$, which are in bijection with the elements of $S_w$. Thus
$h_{K[A_w]}(m) = |S_w(m)|$, and it suffices to show that
$|S_w(m)|$ satisfies (\ref{diffeq}).  We can write
\begin{equation}\label{sred}
S_w(m+d)=\left(S_w(m+d)\, \!\setminus\!\, S_H(m+d)\right)
\dot{\cup} \, S_H(m+d),
\end{equation}
where $S_H(m+d)= \bigcup_{w_i} S_{w_i}(m+d)$, the union being over
the divisors $X(w_i)$ of $X(w)$.  Following the identical
arguments used to deduce (\ref{smh}) (replacing ``$SM$'' by
``$S$'' everywhere), one obtains
\begin{equation}
|S_H(m+d)| =\!\! \sum_{w' < w}\!a_{w'}|S_{w'}(m+d)|,
\end{equation}
for the same integers $a_{w'}$ as in (\ref{smh}).

Establishing an explicit bijection between $S_w(m+d) \!\setminus\!
S_H(m+d)$ and $S_w(m)$ completes the proof, for then (taking
cardinalities of both sides of (\ref{sred})), one sees that
$h_{A_w}(m)$ satisfies (\ref{diffeq}) for all $m \in
\mathbb{Z}_{\geq\, 0}$.

In view of the discussion of flat deformations in Section
\ref{flatdef}, Lemma \ref{main} also implies
\begin{cor}
The set $\{f_{\theta}, \theta \nleq w\} \subset K[O^-]$ forms a
Gr\"obner basis for the ideal it generates.
\end{cor}

\subsection{Combinatorial Interpretation}\label{combin}

We call a multiset $S$ of $(R^+ \!\setminus\! R_{P_d}^+)^{\ast}$ a
{\it t-multipath}, if the chainlength of $S$ is $t$. If $S$ has no
repeated elements (i.e. it is a uniset), then we call it a {\it
t-unipath}. Define $s \in S$ to be a {\it chain-maximal element of
S} if there is no element in $S$ strictly greater than $s$ which
commutes with $s$. Any t-multipath S can be written in the
following manner as the union of $t$ nonintersecting 1-multipaths:
if $S_i$ is the $i^{\textrm{th}}$ 1-multipath, then $S_{i+1}$ is
the multiset of chain-maximal elements (including repetitions) of
$S \!\setminus\! \cup_{k = 1}^i S_k$ (for $i = 0, \cdots, t-1$,
where $S_0$ is defined to be the empty set). If the t-multipath S
is a t-unipath, then each $S_i$ will be a 1-unipath.

Fix $w \in W/W_{P_d}$.  There is a unique expression $w =
s_{\alpha_{i_1}}\!\! \cdots s_{\alpha_{i_{d_w}}}$ such that
$s_{\alpha_{i_k}}\!\! > s_{\alpha_{i_{k+1}}}$ for all k, and all
the reflections pairwise commute; $d_w$ is called the {\it degree}
of $w$.

\begin{ex} Let $w = (3, 5, 7, 8) \in I_{4, 8}$. Then $w = s_{(8,
1)}\,s_{(7,2)}\,s_{(5,4)}$, where $s_{(8, 1)} > s_{(7,2)} > s_{(5,
4)}$ is a chain of commuting reflections.  Thus $d_w = 3$.
\end{ex}

Let $H_j = \{\alpha \in R^+ \!\setminus\! R_{P_d}^+ \, | \,
s_{\alpha} \le s_{\alpha_{i_j}}\}$.  We say that a t-multipath $S$
is w-good if, when written as the union of weighted 1-multipaths
$\cup_{k = 1}^t S_k$ as above, we have that the elements of $S_j$
are in $H_j, j = 1, \cdots, t$. Any multiset in $(R^+
\!\setminus\! R_{P_d}^+)^{\ast}$ is a t-multipath for some t; it
is said to be w-good if the corresponding t-multipath is w-good.

It can be seen that the combinatorial property that a multiset
(resp. uniset) $S$ of $(R^+ \setminus R_{P_d}^+)^{\ast}$ is w-good
is equivalent to the group-theoretic property that $S \in S_w$
(resp. $S \in S_w'$).  Thus Theorem \ref{main1} is equivalent to
the assertion that $h_{\textrm{TC}_{id}X(w)}(m)$ is the number of
w-good multisets of $(R^+ \setminus R_{P_d}^+)^{\ast}$ of degree
m. Letting $M$ be the maximum cardinality of a w-good uniset,
Theorem \ref{main2} is equivalent to the assertion that
$\textrm{mult}_{id}X(w)$ is the number of w-good unisets of
cardinality $M$.

\begin{ex}
Let $w = s_{(15,2)}\,s_{(13,4)}\,s_{(10, 5)} \in  I_{7, 16}$. We
have that $s_{(15,2)} > s_{(13,4)} > s_{(10, 5)}$ is a chain of
commuting reflections, and thus $d_w = 3$.

The diagram below shows the lattice $R^+ \!\setminus\!
R_{P_{\,7}}^+$, where the reflection $s_{(i,j)}$ is denoted by
$i,j$. The set $S$ of reflections which lie along the three
broken-line paths is an example of a $w$-good uniset of maximum
cardinality. In fact, any $w$-good uniset of maximum cardinality
can be seen as the set of reflections lying on three paths in the
lattice, satisfying the following properties:
\begin{itemize}
\item One path starts and ends at ``X'', the second at ``Y'', and
the third at ``Z''.
\item Each path can move only down or to the right.
\item The paths do not intersect.
\end{itemize}
Thus the number of ways of drawing three such paths is
$\textrm{mult}_{id}X(w)$.
\end{ex}


\begin{figure}
\centering
\includegraphics[width=1\textwidth]{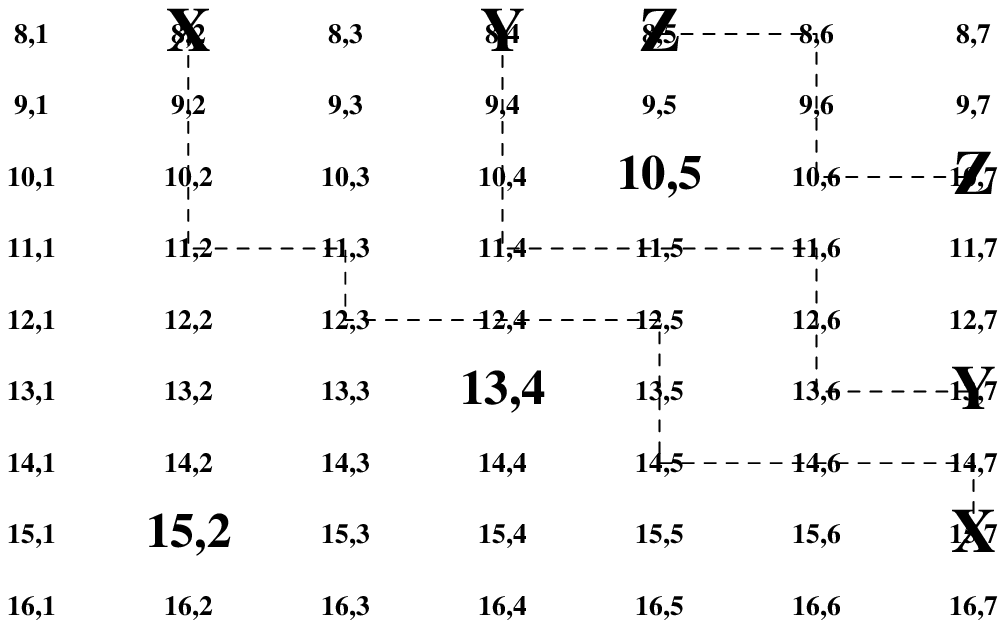}
\caption[]{}
\label{tpaths}
\end{figure}


\section{Conjectures on the Behavior at Other Points}

Let $w, \tau \in W / W_{P_d}$.  Define $S_{w, \tau}$ to be the
multisets $S$ of $(R^+ \setminus R_{P_d}^+)^{\ast}$, such that for
every chain of commuting reflections $s_{\alpha_1} > \cdots >
s_{\alpha_t},\  s_{\alpha_i} \in S$, we have that $w \geq \tau
s_{\alpha_1} \cdots s_{\alpha_t}$; define $S_{w, \tau}'$ to be the
unisets of $(R^+ \setminus R_{P_d}^+)^{\ast}$ having the same
property.  For $m$ a positive integer, define $$S_{w, \tau}(m) =
\{ S\in S_{w, \tau} \, : \, |S|=m\}$$ $$S_{w, \tau}'(m) = \{ S\in
S_{w, \tau}' \, : \, |S|=m\}.$$

We state two conjectures.  First, the Hilbert function
$h_{\textrm{TC}_{\tau}X(w)}(m)$ of the tangent cone to $X(w)$ at
$\tau$ is given by

\begin{conj}
$h_{\textrm{TC}_{\tau}X(w)}(m) = |S_{w, \tau}(m)|, m \in
\mathbb{N}$.
\end{conj}
Second, letting $M$ denote the maximum cardinality of an element
of $S_{w, \tau}'$, the multiplicity $\textrm{mult}_{\tau}X(w)$ of
$X(w)$ at $\tau$ is given by

\begin{conj}
$\textrm{mult}_{\tau}X(w) = |\{S \in S_{w, \tau}' \, : \,
|S|=M\}|$.
\end{conj}


\end{document}